\documentclass[english,onecolumn]{IEEEtran}
\usepackage[T1]{fontenc}
\usepackage[latin9]{inputenc}
\usepackage{geometry}
\usepackage{color}
\usepackage{float}
\usepackage{amsmath}
\usepackage{amsthm}
\usepackage{amssymb}
\usepackage{graphicx}

\makeatletter

\floatstyle{ruled}
\newfloat{algorithm}{tbp}{loa}
\providecommand{\algorithmname}{Algorithm}
\floatname{algorithm}{\protect\algorithmname}

\theoremstyle{plain}
\newtheorem{thm}{\protect\theoremname}
\theoremstyle{definition}
\newtheorem{problem}{\protect\problemname}
\theoremstyle{plain}
\newtheorem{prop}{\protect\propositionname}
\theoremstyle{plain}
\newtheorem{lem}{\protect\lemmaname}

\usepackage{color}

\makeatother

\DeclareMathOperator*{\argmin}{arg\,min}

\usepackage{babel}
\providecommand{\lemmaname}{Lemma}
\providecommand{\problemname}{Problem}
\providecommand{\propositionname}{Proposition}
\providecommand{\theoremname}{Theorem}

\begin{document}
\title{Technical Report: Distributed Asynchronous Large-Scale Mixed-Integer Linear Programming
via Saddle Point Computation}
\author{Luke Fina and Matthew Hale$^{*}$\thanks{$^{*}$Department of Mechanical and Aerospace Engineering, University
of Florida, Gainesville, FL USA. Emails: \texttt{\{l.fina,matthewhale\}@ufl.edu}.
This work was supported by AFOSR under grant FA9550-19-1-0169, ONR
under grants N00014-19-1-2543, N00014-21-1-2495, N00014-22-1-2435,
and by AFRL under grant
FA8651-08-D-0108 TO48.}}
\maketitle
\begin{abstract}
We solve large-scale mixed-integer linear programs (MILPs) via distributed
asynchronous saddle point computation. This is motivated by the MILPs
being able to model problems in multi-agent autonomy, e.g., task assignment problems
and trajectory planning with collision avoidance constraints in multi-robot
systems. To solve a MILP, we relax it with a nonlinear program approximation
whose accuracy tightens as the number of agents increases relative
to the number of coupled constraints. Next, we form an equivalent
Lagrangian saddle point problem, and then regularize 
the Lagrangian in both the
primal and dual spaces to create a regularized Lagrangian that is
strongly-convex-strongly-concave. We then develop a parallelized algorithm
to compute saddle points of the regularized Lagrangian. This algorithm
partitions problems into blocks, which are either scalars or sub-vectors
of the primal or dual decision variables, and it is shown to tolerate
asynchrony in the computations and communications of primal and dual
variables. Suboptimality bounds and convergence rates are presented
for convergence to a saddle point. The suboptimality bound includes
(i) the regularization error induced by regularizing the Lagrangian
and (ii) the suboptimality gap between solutions to the original MILP
and its relaxed form. Simulation results illustrate these theoretical
developments in practice, and show that relaxation and regularization
together have only a mild impact on the quality of solution obtained. 
\end{abstract}

\section{Introduction\label{sec:Introduction}}

Numerous problems in autonomy can be formulated into mixed-integer
programs in which the integer constraints lead to a more realistic
model of the system than a problem without integer constraints. In
particular, mixed-integer programs have been used to model piece-wise
affine functions which can approximate nonlinear dynamical systems
for control \cite{richards2005mixed}, deep neural networks with ReLU
activation functions to improve robustness \cite{fischetti2018deep},
task assignment problems \cite{hendricksondecentralized}, and trajectory
planning problems with collision avoidance \cite{ioan2021mixed}.
In a task assignment problem, the integer decision variables model
discrete decisions about which agents work on which tasks. For collision
avoidance, the integer variables can model a halfspace avoidance region
\cite{ioan2021mixed}. Mixed integer programs are useful in their
modeling power, but scale poorly in comparison to non-discrete optimization
methods due to the NP-Hardness of solving problems with integrality
constraints \cite{conforti2014integer}. Therefore, methods have been
developed to computationally reach an approximate solution. Most existing
methods are centralized, but as the size of problems increases parallelism
and decentralized schemes are needed to partition the problems into
smaller pieces and accelerate computations

To help address this need, in this work we solve large-scale mixed-integer
linear programs (MILPs) via a network of agents in a distributed fashion.
This distributed form allows scaling to large networks of agents in
which each agent performs computations locally and asynchronously shares
the results with other agents. The goal is to minimize a linear cost
while obeying individual constraints, which constrain a single agent,
and shared coupling constraints, which jointly constrain all agents.
This type of problem has been referred to in the literature as ``constraint-coupled
MILPs'' \cite{camisa2021distributed}. For this class of problems,
we seek to find a feasible approximate solution that has bounded,
quantifiable suboptimality using a Lagrangian relaxation approach.

In particular, we solve these large-scale MILPs via distributed asynchronous
computation of saddle points of the Lagrangian. 
First we relax the original
MILP with a nonlinear program approximation whose suboptimality
gap tightens as the number of agents increases relative to the number
of coupled constraints. Then we form a Lagrangian saddle point problem
equivalent to that nonlinear program and regularize in the primal
and dual spaces to make the Lagrangian strongly-convex-strongly-concave.
Regularization not only provides robustness to asynchronous updates,
but also ensures the uniqueness
of a solution; this uniqueness is not guaranteed in unregularized cases, unless
it is assumed outright. 

We develop a parallelized saddle point finding algorithm
and apply it to the regularized Lagrangian. This algorithm
partitions problems into blocks that are either scalars or sub-vectors
of the primal or dual decision variables, and this algorithm is shown
to be tolerant to asynchrony in computations and communications. The
block structure simplifies the needed communications for our large-scale
MILP problems because these problems are separable, which implies
that agents can compute their local gradients without needing to share
information with each other. We present convergence analysis for our
algorithm, along with a suboptimality bound between the solution to
the original MILP and the approximate solution obtained using our
algorithm. 

Our particular algorithm draws on recent developments in distributed
primal-dual optimization and solving large-scale mixed-integer linear
programs. Similar work appears in \cite{camisa2021distributed,camisa2019primal},
but one major distinction is that those works use a consensus-based
method while we use a block-based method. Consensus methods require
graph connectivity of some form over time, which typically cannot
be verified locally by the agents in the network. Block-based methods,
on the other hand, require only bounded delays along each communication
channel, which are typically easier to verify. Block-based optimization
can also lead to the communication of less information between agents
and fewer required computations per agent, since each agent computes
updates to and communicates only a (usually small) subset of the decision
variables in a problem. This is in contrast to consensus-type methods,
in which the entire decision vector is updated and communicated by
each agent. These differences are particularly advantageous for multi-robot
systems, especially in environments with unreliable communications,
and in bandwidth-limited systems, such as the large-scale systems
we consider.

Distributed convex constrained primal-dual optimization is a foundational
element for the current work because it is used to solve the relaxed
MILP \cite{hendrickson2021totally}. Earlier works on decentralized
primal-dual constrained optimization address the effect of delay bounds,
but \cite{hale2017asynchronous} did not address integer programming
and \cite{falsone2017dual} used consensus-based methods. In multi-agent
optimization, few works solve large-scale mixed integer linear programs,
but of those that do most recent works are based on the developments
in \cite{vujanic2016decomposition}. That work was later improved
for a tighter suboptimality bound in \cite{falsone2017dual} and then
it was subsequently distributed in \cite{falsone2020distributed}.
Our algorithm differs from that in \cite{falsone2020distributed}
because our algorithm is block-based and asynchronous while the algorithm
in \cite{falsone2020distributed} is consensus-based and synchronous.
Our algorithm converges with bounded delays, and the rate of convergence
depends on the length of the delays and the regularization terms in
a way that we quantify.

To summarize, the main contributions of this paper are the following:
\begin{itemize}
\item We are the first to solve large-scale MILPs with block-based asynchronous
communications in a decentralized algorithm (Algorithm \ref{alg:algorithm1}).
\item We are the first to bound the suboptimality gap between the original
MILP and its approximation by a regularized Lagrangian saddle point formulation (Theorem
\ref{thm:SuboptimalityBound}). 
\item We regularize to remove the unique linear program assumption present
in all other large-scale MILP works.
\item We present a rate of convergence for our algorithm under partial asynchrony
(Theorem \ref{thm:overallConvergence}).
\end{itemize}
This paper is organized as follows. Preliminaries and problem statements
are given in Section \ref{sec:Preliminaries-and-Problem}. Section
\ref{sec:Distributed-Algorithm} presents the algorithm, convergence
analysis, and suboptimality bounds. Section \ref{sec:Simulation}
demonstrates the algorithm in simulation, and Section \ref{sec:Conclusion}
concludes.

\section{Preliminaries and Problem Statements\label{sec:Preliminaries-and-Problem}}

This section establishes the problems we solve, the assumptions on
those problems, and the relationship between the relaxed and unrelaxed
forms of the MILPs we solve. Below, we use the notation $conv(\cdot)$
to denote the convex hull of a set, the set of indices $[s]\triangleq\{1,...,s\}$
for any $s\in\mathbb{N}$, $\mathbb{R}$ to denote real numbers, $\mathbb{Z}$
to denote integers, $\mathbb{R}_{+}$ to denote positive real numbers,
$B_o(r)$ to denote the closed ball of radius~$r$ centered on the origin, 
and $\Pi_{S}[\cdot]$ to denote the Euclidean projection onto a compact,
convex set $S$. 

We study the following MILP over a network of $m$ agents.
\begin{problem}
\emph{(Original MILP)\label{prob:Original-MILP}}
\end{problem}
\begin{align*}
\underset{x_{1},...,x_{m}}{\text{minimize}} & \enskip f(x)=\sum_{\ell=1}^{m}c_{\ell}^{T}x_{\ell}\\
\text{subject to} & \enskip g(x)=\sum_{\ell=1}^{m}A_{\ell}x_{\ell}\leq b\\
 & x_{\ell}\in X_{\ell}\text{\ensuremath{\quad}for all }\ell\in[m],
\end{align*}
\emph{where $X_{\ell}\triangleq\{x_{\ell}\in\mathbb{Z}^{p_{\ell}}\times\mathbb{R}^{q_{\ell}}:S_{\ell}x_{\ell}\leq s_{\ell}\}$
for a given $S_{\ell}\in\mathbb{R}^{v\times(p_{\ell}+q_{\ell})}$
and }$s_{\ell}\in\mathbb{R}^{v}$\emph{, the linear cost is given
by $c_{\ell}\in\mathbb{R}^{p_{\ell}+q_{\ell}}$, the coupling constraints
are given by $A_{\ell}\in\mathbb{R}^{y\times(p_{\ell}+q_{\ell})}$,
and $b\in\mathbb{R}^{y}$. }

To solve Problem \ref{prob:Original-MILP}, branch-and-bound or cutting
plane methods could be used to achieve an exact optimal solution.
However, these methods are computationally expensive and can be too
slow for many real world applications, especially in large-scale instances.
To this end, the following problem approximates Problem \ref{prob:Original-MILP}
and we take advantage of the separable problem structure in a later
proposition. Specifically, we obtain an approximate form of Problem
\ref{prob:Original-MILP} where $\rho$ is a vector that tightens
the coupling constraints in Problem \ref{prob:Relaxed-MILP}, which
will assist in reconstructing
a feasible solution to Problem \ref{prob:Original-MILP}. 
\begin{problem}
\emph{(Relaxed MILP; approximation of Problem \ref{prob:Original-MILP})\label{prob:Relaxed-MILP}}
\end{problem}
\textcolor{black}{
\begin{align*}
\underset{z}{\text{minimize}} & \enskip f(z)=c^{T}z\\
\text{subject to} & \enskip g(z)=Az\leq b-\rho\\
 & \enskip z\in Z,
\end{align*}
}\textcolor{black}{\emph{where $z=[z_{1}^{T},...,z_{m}^{T}]^{T}\in\mathbb{R}^{p_{1}+q_{1}+\cdots+p_{m}+q_{m}}$,
$Z=Z_{1}\times\cdots\times Z_{m}\subseteq\mathbb{R}^{p_{1}+q_{1}+...+p_{m}+q_{m}}$,
$Z_{\ell}=conv(X_{\ell}),z_{\ell}\in\mathbb{R}^{p_{\ell}+q_{\ell}},$
$A=[A_{1},A_{2},...,A_{m}]\in\mathbb{R}^{y\times(p_{1}+q_{1}+...+p_{m}+q_{m})}$,
$A_{\ell}\in\mathbb{R}^{y\times(p_{\ell}+q_{\ell})}$ $,c=[c_{1}^{T},...,c_{m}^{T}]^{T}\in\mathbb{R}^{p_{1}+q_{1}+...+p_{m}+q_{m}},$
$b\in\mathbb{R}^{y}$, and $\rho\in\mathbb{R}_{+}^{y}$.}}

\noindent The symbol $z$ is introduced to differentiate between the
decision variables in Problem \ref{prob:Original-MILP} and those in Problem
\ref{prob:Relaxed-MILP}.\emph{ }We are interested in large-scale
problems, which are those with $m\gg y$. We make the following assumptions
about Problem \ref{prob:Relaxed-MILP}.

\noindent \textbf{Assumption 1 (Nonempty and Compact) }\emph{The set
$Z$ is nonempty and compact.}

\noindent This assumption ensures the existence of a solution and
allows for distributed projected update laws which will ensure satisfaction
of the agents' set constraints.

\noindent \textbf{Assumption 2 (Slater's Condition)} \emph{Slater's
condition holds for $g,$} i.e., there exists $\bar{z}\in Z$ such
that $g(\bar{z})<0.$

\noindent This assumption will simplify our convergence analysis by
allowing us to formulate a saddle point problem over a compact domain. 

We solve Problem \ref{prob:Relaxed-MILP} using a primal-dual approach.
Specifically, we will use a decentralized form of the classic Uzawa
iteration, also called gradient descent-ascent, to find a saddle point
of the Lagrangian associated with Problem \ref{prob:Relaxed-MILP}.
This Lagrangian is $L(z,\lambda)=c^{T}z+\lambda^{T}(Az-b+\rho)$,
which is affine in both $x$ and $\lambda.$ Thus it is convex in
$x$ and concave in $\lambda.$ To ensure convergence of a decentralized
asynchronous implementation of the Uzawa iteration, we will regularize
$L$ to make it strongly convex in $x$ and strongly concave in $\lambda$.
These properties have both been shown to improve convergence of decentralized
algorithms \cite{koshal2011multiuser}. We apply a Tikhonov regularization
in the primal and dual spaces to find

\textcolor{black}{
\begin{align}
 & L_{\kappa}(z,\lambda)=c^{T}z+\frac{\alpha}{2}\vert\vert z\vert\vert^{2}+\lambda^{T}(Az-b+\rho)-\frac{\delta}{2}\vert\vert\lambda\vert\vert^{2},\label{eq:RegularizedLagran}
\end{align}
}where $\lambda\in\mathbb{R}_{+}^{y}$ is the dual variable, $\alpha>0$
is the primal regularization term, $\delta>0$ is the dual regularization
term, and we use $\kappa=(\alpha,\delta$) to simplify notation. 

Strong duality holds for Problem \ref{prob:Relaxed-MILP} \cite[Section 5.2.3]{boyd2004convex},
which means it can be solved by computing the saddle point of the
unregularized Lagrangian $L$. Instead of doing so exactly, we will
compute a saddle point of the regularized Lagrangian $L_{\kappa}$
with $\alpha>0$ and $\delta>0$ set to small values. 
As we show in Section IV, this regularization
induces only minor error in the solutions we obtain. \textcolor{black}{Another
advantage of regularizing is that we do not need to assume uniqueness
of the solution to Problem \ref{prob:Relaxed-MILP} since we have
strong convexity from regularization and thus $L_{\kappa}$ has a
unique saddle point. This is unlike existing large-scale MILP methods
that rely on assuming that there is a unique solution to solve Problem~\ref{prob:Relaxed-MILP}
\cite{camisa2021distributed,falsone2017dual,vujanic2016decomposition}. }

\textcolor{black}{Formally, we will solve Problem }\ref{prob:Relaxed-MILP}
by finding a saddle point of $L_{\kappa}$, namely

\[
(\hat{z}_{\kappa},\hat{\lambda}_{\kappa})=\underset{z\in Z}{\text{arg min}}\underset{\lambda\in\mathbb{R}_{+}^{y}}{\enskip\text{arg max}}\enskip c^{T}z+\frac{\alpha}{2}\vert\vert z\vert\vert^{2}+\lambda^{T}(Az-b+\rho)-\frac{\delta}{2}\vert\vert\lambda\vert\vert^{2}.
\]

Bounding the suboptimality gap between
Problem \ref{prob:Original-MILP} and the saddle point $(\hat{z}_{\kappa},\hat{\lambda}_{\kappa}$)
is needed to show that solving the relaxed saddle point problem will
approximately solve the original MILP. To that end, we first upper
bound the dual variable, which will later allow us to bound the suboptimality
that results from the relaxation and provide an upper bound on the maximum
norm of the dual variable. We now explicitly define $\rho$ as the
contraction vector defined by \cite{vujanic2016decomposition} as

\begin{align}
\rho_{\zeta}\triangleq y\cdot\underset{\ell\in[m]}{\max}\left(\underset{x_{\ell}\in X_{\ell}}{\max}A_{row,\ell}x_{\ell}-\underset{x_{\ell}\in X_{\ell}}{\min}A_{row,\ell}x_{\ell}\right),\label{eq:RhoEqn}
\end{align}
where $\ensuremath{\rho\in\mathbb{R}_{+}^{y}}$, $A_{row,\ell}$ is
the $\ell^{th}$ row of A, and $\rho_{\zeta} > 0$ is the $\zeta^{th}$
entry of $\rho$ for~$\zeta \in [y]$. 
We next define a radius
$r$ to satisfy $Z\subseteq B_{o}(r)$. This value of $r$ is a bound
on the maximum norm of $z\in Z,$ and it can be found via
\begin{align}
r\triangleq\underset{z\in Z}{\max}\enskip & \vert\vert z\vert\vert.\label{eq:feasProblem-1}
\end{align}
We will use the following proposition for the subsequent suboptimality
theorem. 
\begin{prop}
\label{prop:relationZandXTheorem} Let Problem \ref{prob:Relaxed-MILP}
have a feasible solution and let $\tilde{z}$ be a vertex of its feasible
set. Then, there exists an index set $I_{\mathbb{Z}}\subseteq [m]$
with cardinality $\vert I_{\mathbb{Z}}\vert\geq m-y$ such
that $\tilde{z}_{\ell}\in X_{\ell}$ for all $\ell\in I_{\mathbb{Z}}.$ 
\end{prop}
\begin{IEEEproof}
See \cite[Theorem 3.3]{vujanic2016decomposition}.
\end{IEEEproof}
This proposition means the optimal solution to Problem \ref{prob:Relaxed-MILP}
is partially mixed-integer and there is a bound on the number of agents
whose solution may not be mixed-integers for Problem \ref{prob:Original-MILP}.
In particular, as the size of problems grows such that $m\gg y$,
we find that the integrality constraints become approximated better.
In the following lemma, we bound the dual multipliers to be used in
a later theorem.
\begin{lem}
\label{lem:saddlePoint}Let Assumptions 1 and 2 hold. Then

\[
\hat{\lambda}_{\kappa}\in M\triangleq\left\{ \lambda\in\mathbb{R}_{+}^{y}:\vert\vert\lambda\vert\vert_{1}\leq\frac{c^{T}\bar{z}+\frac{\alpha}{2}\vert\vert\bar{z}\vert\vert^{2}-\vert\vert c\vert\vert\cdot r}{\underset{1\leq j\leq m}{\min}-A_{row,j}\bar{z}+\rho_{j}-b_{j}}\right\} ,
\]
where $\bar{z}$ is a Slater point of $g$, \textup{$c^{T}\bar{z}$}
is the cost at the Slater point, $\alpha$ is the primal regularization
term, $r$ is the radius found from 
(\ref{eq:feasProblem-1}), $A_{row,j}$ denotes row $j$ of $A$,
and $-A_{j}\bar{z}+\rho_{j}-b_{j}$ is the $jth$ entry of $g(\bar{z})$.
\end{lem}
\begin{IEEEproof}
As discussed in \cite[Section II.C]{hale2015cloud} for any Slater
point $\bar{z}$ of $g$ we have

\begin{equation}
\hat{\lambda}_{\kappa}\in M\triangleq\left\{ \lambda\in\mathbb{R}_{+}^{y}:\vert\vert\lambda\vert\vert_{1}\leq\frac{f(\bar{z})+\frac{\alpha}{2}\vert\vert\bar{z}\vert\vert^{2}-\underset{z\in Z}{\min}f(z)}{\underset{1\leq j\leq m}{\min}-g_{j}(\bar{z})}\right\} .\label{eq:Lem4_part1}
\end{equation}
We expand $g(\bar{z})=A\bar{z}+\rho-b$, which implies $g_{j}(\bar{z})=A_{j}\bar{z}+\rho_{j}-b_{j}.$
Then 
\begin{equation}
\underset{1\leq j\leq m}{\min}-g_{j}(\bar{z})=\underset{1\leq j\leq m}{\min}-A_{row,j}\bar{z}-\rho_{j}+b_{j}.\label{eq:Lem4_part2}
\end{equation}
Next, we find a lower bound for $\underset{z\in Z}{\min}f(z)$. It
follows that

\begin{equation}
\underset{z\in Z}{\min}\quad c^{T}z\geq\underset{z\in B_{o}(r)}{\min}c^{T}z=c^{T}\cdot\begin{pmatrix}-\frac{c}{\vert\vert c\vert\vert}\end{pmatrix}\cdot r=-\vert\vert c\vert\vert\cdot r.\label{eq:Lem4_part3}
\end{equation}
Substituting (\ref{eq:Lem4_part2}) and (\ref{eq:Lem4_part3})
into (\ref{eq:Lem4_part1}) gives the desired bound.
\end{IEEEproof}
We now bound the suboptimality in the primal component of the saddle
point of $L_{\kappa}$ that is due to the regularization of $L.$
Below, we will use this lemma to bound the suboptimality of the primal
component of the saddle point of $L_{\kappa}$, namely $\hat{z}_{\kappa}$,
relative to the optimal solution to Problem \ref{prob:Original-MILP}.
\begin{lem}
\label{lem:upperboundLagranCost-1} Let Assumptions 1 and 2 hold.
Then

\[
c^{T}\hat{z}_{\kappa}\leq c^{T}z^{*}+\vert\vert c\vert\vert\cdot\left(\frac{c^{T}\bar{z}+\frac{\alpha}{2}\vert\vert\bar{z}\vert\vert^{2}-\vert\vert c\vert\vert\cdot r}{\underset{1\leq j\leq m}{\min}-A_{row,j}\bar{z}+\rho_{j}-b_{j}}\right)\cdot\sqrt{\frac{\delta}{2\alpha}}+\frac{\alpha}{2}\cdot r,
\]
where $\hat{z}_{\kappa}$ is the primal component of the saddle point
of the regularized Lagrangian, $z^{*}$ is the primal component of
the saddle point of the unregularized Lagrangian~$L$,
$\bar{z}$ is a Slater point for $g$, $\alpha$ is the primal regularization
term, $\delta$ is the dual regularization term, $r$ is found from
 (\ref{eq:feasProblem-1}), and \textup{$-A_{j}\bar{z}-\rho_{j}+b_{j}$
is the explicit form of $-g_{j}(\bar{z})$.}
\end{lem}
\begin{IEEEproof}
From the general form presented in \cite[Lemma 3.3]{koshal2011multiuser},
we find an upper bound between the costs at $\hat{z}_{\kappa}$ and
$z^{*}$ via
\begin{equation}
\vert f(\hat{z}_{\kappa})-f(z^{*})\vert\leq\underset{z\in Z}{\max}\vert\vert\nabla f(z)\vert\vert\cdot\underset{\lambda\in M}{\max}\vert\vert\lambda\vert\vert\cdot\sqrt{\frac{\delta}{2\alpha}}+\frac{\alpha}{2}\cdot\underset{z\in Z}{\max}\vert\vert z\vert\vert.\label{eq:optBoundFull}
\end{equation}
From Problem \ref{prob:Relaxed-MILP}, we find
\[
\underset{z\in Z}{\max}\vert\vert\nabla f(z)\vert\vert=\vert\vert c\vert\vert,
\]
since $\nabla f \equiv c.$ Then we find an upper bound for $\underset{z\in Z}{\max}\vert\vert z\vert\vert$
as

\begin{equation}
\underset{z\in Z}{\max}\vert\vert z\vert\vert\leq\underset{z\in B_{o}(r)}{\max}\vert\vert z\vert\vert=r,\label{eq:maxR}
\end{equation}
which can be found explicitly from 
(\ref{eq:feasProblem-1}). Since $\ensuremath{\lambda}$ is a vector,
we always have $\underset{\lambda\in M}{\max}\vert\vert\lambda\vert\vert_{2}\leq\underset{\lambda\in M}{\max}\vert\vert\lambda\vert\vert_{1}.$
From Lemma 1, we can therefore upper bound $\vert\vert\lambda\vert\vert_{2}$
as
\begin{align}
 \underset{\lambda\in M}{\max}\vert\vert\lambda\vert\vert_{2} &\leq\underset{\lambda\in M}{\max}\left(\frac{c^{T}\bar{z}+\frac{\alpha}{2}\vert\vert\bar{z}\vert\vert^{2}-\vert\vert c\vert\vert\cdot r}{\underset{1\leq j\leq m}{\min}-A_{row,j}\bar{z}+\rho_{j}-b_{j}}\right) \\ &=\frac{c^{T}\bar{z}+\frac{\alpha}{2}\vert\vert\bar{z}\vert\vert^{2}-\vert\vert c\vert\vert\cdot r}{\underset{1\leq j\leq m}{\min}-A_{row,j}\bar{z}+\rho_{j}-b_{j}}.\label{eq:boundmaxLamb-1}
\end{align}
Then, using (\ref{eq:maxR}) and (\ref{eq:boundmaxLamb-1}) in (\ref{eq:optBoundFull}) we find
\begin{equation}
\vert c^{T}\hat{z}_{\kappa}-c^{T}z^{*}\vert\leq\vert\vert c\vert\vert\cdot\left(\frac{c^{T}\bar{z}+\frac{\alpha}{2}\vert\vert\bar{z}\vert\vert^{2}-\vert\vert c\vert\vert\cdot r}{\underset{1\leq j\leq m}{\min}-A_{row,j}\bar{z}+\rho_{j}-b_{j}}\right)\cdot\sqrt{\frac{\delta}{2\alpha}}+\frac{\alpha}{2}\cdot r.
\end{equation}
Rearranging terms we have an upper bound on $c^{T}\hat{z}_{\kappa}$
via 

\[
c^{T}\hat{z}_{\kappa}\leq c^{T}z^{*}+\vert\vert c\vert\vert\cdot\left(\frac{c^{T}\bar{z}+\frac{\alpha}{2}\vert\vert\bar{z}\vert\vert^{2}-\vert\vert c\vert\vert\cdot r}{\underset{1\leq j\leq m}{\min}-A_{row,j}\bar{z}+\rho_{j}-b_{j}}\right)\cdot\sqrt{\frac{\delta}{2\alpha}}+\frac{\alpha}{2}\cdot r.
\]
\end{IEEEproof}
Now that we have Lemma \ref{lem:saddlePoint}, Lemma \ref{lem:upperboundLagranCost-1},
and Proposition \ref{prop:relationZandXTheorem}, we construct the
suboptimality bound on the difference in optimal costs between Problem
\ref{prob:Original-MILP} and the primal component
of the regularized Lagrangian,~$\hat{z}_{\kappa}$. 
\begin{thm}
\label{thm:SuboptimalityBound} (Suboptimality Bound) Let Assumptions
1 and 2 hold. Let\textup{\emph{ $x^{*}$ denote the solution to Problem
\ref{prob:Original-MILP}, and let }}\textup{$\hat{z}_{\kappa}$}
denote the primal component of the saddle point of~$L_{\kappa}$ with~$\rho$
as defined in (\ref{eq:RhoEqn}). Let \emph{$\eta_{\ell}=\underset{x_{\ell}\in X_{\ell}}{\max}c_{\ell}^{T}x_{\ell}-\underset{x_{\ell}\in X_{\ell}}{\min}c_{\ell}^{T}x_{\ell}$}
and let~$\bar{z} \in Z$ denote a Slater point for~$g$. 
Then 
\[
c^{T}\hat{z}_{\kappa}-c^{T}x^{*}\leq y\cdot\underset{\ell\in[m]}{\max}\enskip\eta_{\ell}+\vert\vert c\vert\vert\cdot\left(\frac{c^{T}\bar{z}+\frac{\alpha}{2}\vert\vert\bar{z}\vert\vert^{2}-\vert\vert c\vert\vert\cdot r}{\underset{1\leq j\leq m}{\min}-A_{row,j}\bar{z}+\rho_{j}-b_{j}}\right)\cdot\sqrt{\frac{\delta}{2\alpha}}+\frac{\alpha}{2}\cdot r,
\]
where~$\alpha > 0$ is the primal regularization term,~$\delta > 0$
is the dual regularization term, and~$r$ is from (\ref{eq:feasProblem-1}). 
\end{thm}
\begin{IEEEproof}
Let $z^{*}$ denote the solution to Problem \ref{prob:Relaxed-MILP}
with $\rho$ as defined in (\ref{eq:RhoEqn}). 
Expanding $c^{T}\hat{z}_{\kappa}-c^{T}x^{*}$, we have
\begin{align}
c^{T}\hat{z}_{\kappa}-c^{T}x^{*} & =(c^{T}\hat{z}_{\kappa}-c^{T}z^{*})+(c^{T}z^{*}-c^{T}x^{*}).\label{eq:FullSubOptEquation}
\end{align}
According
to Proposition \ref{prop:relationZandXTheorem}, there exists $I_{\mathbb{Z}}$
such that $\vert I_{\mathbb{Z}}\vert\geq m-y$ and $z_{\ell}^{*}\in X_{\ell}$
for each $\ell\in I_{\mathbb{Z}}$. This means $z_{\ell}^{*}=x_{\ell}^{*}$
for all $\ell\in I_{\mathbb{Z}}$. Define $I_{\mathbb{R}}=[m]\backslash I_{\mathbb{Z}}$,
which contains indices such that $z_{\ell}^{*}\in conv(X_{\ell})\backslash X_{\ell}$.
We can simplify the difference between $c^{T}z^{*}$ and $c^{T}x^{*}$ as
\begin{align}
c^{T}z^{*}-c^{T}x^{*} & =\sum_{i\in I_{\mathbb{Z}}}(c_{i}^{T}z_{i}^{*}-c_{i}^{T}x_{i}^{*})+\sum_{j\in I_{\mathbb{R}}}(c_{j}^{T}z_{j}^{*}-c_{j}^{T}x_{j}^{*})\nonumber \\
 & =\sum_{j\in I_{\mathbb{R}}}(c_{j}^{T}z_{j}^{*}-c_{j}^{T}x_{j}^{*}).\label{eq:IRTerms-1}
\end{align}
We bound $c^{T}z^{*}-c^{T}x^{*}$ in terms of the remaining $I_{\mathbb{R}}$
terms by upper and lower bounding each term in (\ref{eq:IRTerms-1}).
Because Problem \ref{prob:Relaxed-MILP} is a relaxation of Problem
\ref{prob:Original-MILP}, we have $c_{\ell}^{T}z_{\ell}^{*}\leq c_{\ell}^{T}x_{\ell}$
for every $x_{\ell}\in X_{\ell}$. Then, for any $x_{\ell}\in X_{\ell},$ we
have $c_{j}^{T}z_{j}^{*}\leq c_{\ell}^{T}x_{\ell}\leq\underset{x_{\ell}\in X_{\ell}}{\max}c_{\ell}^{T}x_{\ell}$.
Next, for every $j\in I_{\mathbb{R}},$ for $x^{*}$ we have $\underset{x_{\ell}\in X_{\ell}}{\min}c_{\ell}^{T}x_{\ell}\leq c_{j}^{T}x_{j}^{*}$.
Then we find 

\[
c_{j}^{T}z_{j}^{*}-c_{j}^{T}x_{j}^{*}\leq\underset{x_{\ell}\in X_{\ell}}{\max}c_{\ell}^{T}x_{\ell}-\underset{x_{\ell}\in X_{\ell}}{\min}c_{\ell}^{T}x_{\ell}.
\]
We thus find an upper bound on the distance between $c^{T}z^{*}$and
$c^{T}x^{*}$ as

\begin{align}
 & \sum_{j\in I_{\mathbb{R}}}c_{j}^{T}z_{j}^{*}-c_{j}^{T}x_{j}^{*}\leq\sum_{j\in I_{\mathbb{R}}}\begin{pmatrix}\underset{x_{\ell}\in X_{\ell}}{\max}c_{\ell}^{T}x_{\ell}-\underset{x_{\ell}\in X_{\ell}}{\min}c_{\ell}^{T}x_{\ell}\end{pmatrix}=\vert I_{\mathbb{R}}\vert\cdot\underset{\ell\in[m]}{\max}\begin{pmatrix}\underset{x_{\ell}\in X_{\ell}}{\max}c_{\ell}^{T}x_{\ell}-\underset{x_{\ell}\in X_{\ell}}{\min}c_{\ell}^{T}x_{\ell}\end{pmatrix}\leq y\cdot\underset{\ell\in[m]}{\max}\enskip\eta_{\ell},\label{eq:thm1Part1Result}
\end{align}
which follows from the definition of $\eta_{\ell}$ and the fact that
$\vert I_{\mathbb{R}}\vert\leq y.$ Next, we use Lemma \ref{lem:upperboundLagranCost-1}
to relate $c^{T}z^{*}$ to $c^{T}\hat{z}_{\kappa}$ via

\begin{equation}
c^{T}\hat{z}_{\kappa}-c^{T}z^{*}\leq\vert\vert c\vert\vert\cdot\left(\frac{c^{T}\bar{z}+\frac{\alpha}{2}\vert\vert\bar{z}\vert\vert^{2}-\vert\vert c\vert\vert\cdot r}{\underset{1\leq j\leq m}{\min}-A_{row,j}\bar{z}+\rho_{j}-b_{j}}\right)\cdot\sqrt{\frac{\delta}{2\alpha}}+\frac{\alpha}{2}\cdot r.\label{eq:thm1Part2Result}
\end{equation}
Using the results of (\ref{eq:IRTerms-1}), (\ref{eq:thm1Part1Result})
and (\ref{eq:thm1Part2Result}) in (\ref{eq:FullSubOptEquation})
completes the proof.
\end{IEEEproof}
Now that we have a suboptimality bound between 
the solution to Problem \ref{prob:Original-MILP}
and the primal component~$\hat{z}_{\kappa}$ of the saddle point of the regularized
Lagrangian~$L_{\kappa}$, 
we will find the saddle point of $L_{\kappa}$ in a parallelized way, 
and the next section develops the algorithm for doing so.

\section{Distributed Algorithm\label{sec:Distributed-Algorithm}}
In this section we define the main algorithm we use to
find the saddle point of~$L_{\kappa}$, 
analyze its convergence, and describe how its output can
be mapped back to an approximate solution to Problem~\ref{prob:Original-MILP}. 

\subsection{Algorithm Definition}
The distributed algorithm we develop consists of three types of operations;
(i) a primal variable update, (ii) a dual variable update, and (iii)
communication of these updates. 
Recall from Problem~\ref{prob:Relaxed-MILP} that~$z = (z_1^T, \ldots, z_m^T)^T$, which 
is partitioned into~$m$ blocks. 
We assume that there are~$H > m$ agents\footnote{One can use~$H = m$ agents and simply have
each agent perform the duties of both a primal agent and a dual agent, though here
we assume that~$H > m$ to allow primal and dual agents to be separate, which
simplifies our discussion.} 
to find the saddle point~$(\hat{z}_{\kappa}, \hat{\lambda}_{\kappa})$, 
of which there are~$m$ primal agents indexed over the set $A_{P}=[m_{1}]$
and $H-m$ dual agents, indexed over the set $A_{D}=[H]\backslash[m]$.
Each primal and dual agent has partial knowledge of
Problem \ref{prob:Relaxed-MILP} and only updates its assigned block
of the decision variables. We partition the dual 
variable~$\lambda \in \mathbb{R}^y_{+}$ into blocks via
\begin{equation}
\lambda = (\lambda_1^T, \ldots, \lambda_{H-m}^T)^T,
\end{equation}
where~$\lambda_c \in \mathbb{R}^{r_c}_{+}$
for~$c \in A_D$ and~$\sum_{c \in A_D} r_c = y$.

Computations and communications are asynchronous, which leads to disagreements
among the values of the decision variables that agents store onboard.
Therefore, there is a need to 
compute the saddle point~$(\hat{z}_{\kappa}, \hat{\lambda}_{\kappa})$
with an update law that is robust to asynchrony. We use a decentralized
form of the classic Uzawa iteration, also called gradient descent-ascent,
to find a saddle point of the regularized Lagrangian~$L_{\kappa}$ associated with
Problem \ref{prob:Relaxed-MILP}. 

Each primal agent updates a block of the primal variables it stores onboard. 
Primal agent $i$ stores onboard
itself a local copy of the vector of primal decision variables, 
denoted $z^{i} \in Z$,
and local copy of the vector of dual decision variables, 
denoted $\lambda^{i} \in M$. 
Dual agent $q$ stores a local
primal vector onboard, denoted $z^{q} \in Z$,
and local dual vector, denoted $\lambda^{q} \in M$. 
Each dual agent updates
its block of the dual variable in the copy it stores onboard.
Each dual variable corresponds to an inequality constraint on the primal
variables, and each dual agent sends values
of its updated block to all primal
agents whose decision variables appear in the constraints that correspond
to those dual variables.

We define $D\subseteq\mathbb{N}$ as the set of times when all dual agents
compute updates to their decision variables, 
and we define~$K^{i}\subseteq\mathbb{N}$ as the set of times
at which primal agent $i\in A_{P}$ computes updates
to its decision variables. To state an algorithm, we use
$k$ as the iteration count used by all the primal agents,
and $t$ as the iteration count used
by all dual agents. 
The sets $K^{i}$ and $D$ are tools for analysis and discussion, and they need not be known by any agent. 
Additionally,
within the vector~$z^i$, the
notation $z_{[j]}^{i}$ indicates agent $i$'s value for the primal
block $j$, $z_{[i]}^{i}$ indicates agent $i$'s onboard value for
the primal block $i$, 
$\lambda_{[j]}^{i}$ indicates agent~$i$'s onboard value for the dual block~$j$, 
$\nabla_{z_{[i]}} := \frac{\partial}{\partial z_{[i]}}$ 
is the derivative with
respect to the $i^{th}$ block of $z$, and
$\nabla_{\lambda_{[q]}} := \frac{\partial}{\partial \lambda_{[q]}}$
is the derivative with respect to the~$q^{th}$ block of $\lambda$.
Then $\lambda^{i}(k)$ denotes the value of~$\lambda$ stored onboard
agent~$i$ at time~$k$, and~$z^i(k)$ denotes the value of~$z$
stored onboard agent~$i$ at time~$k$. 

We make the following assumption about delays, 
termed ``partial asynchrony''~\cite{bertsekas2015parallel}. 

\textbf{Assumption 3 (Bounded Delays) }

\emph{Let $K^{i}$ be the set of times when primal agent i performs
updates and let $\tau_{j}^{i}(k)$ be the time when primal agent j
originally computed the value of $z_{j}^{i}(k)$ onboard by agent
i at time k. Then there exists a positive integer B such that }
\begin{enumerate}
\item \emph{For every $i\in A_{P}$, at least one of the elements of the
set $\{k,k+1,...,k+B-1\}$ is in $K^{i}.$}
\item \emph{There holds $k-B<\tau_{j}^{i}(k)\leq k$, for all $i,j\in [m]$,
$j\neq i,$ and all $k\in K^{i}.$}
\end{enumerate}
This assumption ensures that no primal block onboard
any agent was computed more than~$B$ 
iterations prior to the current time, and it ensures that
all primal agents perform at least one computation every~$B$ timesteps. 

Primal agent $i$ updates its block~$z^i_{[i]}$ via projected gradient descent
at each time $k\in K^{i}$. If $k\notin K^{i}$,
then agent $i$ does not update and~$z^i_{[i]}$ is held constant. 
If a communication of an updated block
from primal agent $j$ is received, then agent $i$ stores it in the block 
$z_{[j]}^{i}$ by overwriting the previous value of~$z^i_{[j]}$ that it had
onboard. 
At times at which agent~$i$ does not receive a communication from primal agent~$j$, 
$z_{[j]}^{i}$ is held constant. At each $t\in D$, dual
agents update via projected gradient ascent onto the set $M_{q}$, which
we define as
\begin{equation}
M_q := \left\{ \lambda_c \in \mathbb{R}_{+}^{r_c}:\vert\vert\lambda_c\vert\vert_{1}\leq\frac{c^{T}\bar{z}+\frac{\alpha}{2}\vert\vert\bar{z}\vert\vert^{2}-\vert\vert c\vert\vert\cdot r}{\underset{1\leq j\leq m}{\min}-A_{row,j}\bar{z}+\rho_{j}-b_{j}}\right\}. 
\end{equation}
We note that the upper bound in~$M_q$ is the same as that in~$M$ in Lemma~\ref{lem:saddlePoint} and
thus is a bound on each individual block of~$\lambda$ as well. 

Primal agents $i$'s computations take the form
\begin{gather*}
z_{[i]}^{i}(k+1)=\begin{cases}
\Pi_{Z_{i}}[z_{[i]}^{i}(k)-\gamma\nabla_{z_{[i]}}L_{\kappa}(z^{i}(k),\lambda^{i}(k))] & k\in K^{i}\\
z_{i}^{i}(k) & \text{otherwise.}
\end{cases}
\end{gather*}
For all~$q \in A_D$, dual agent~$q$ only performs an update
occasionally. Specifically, it only performs an update after
the counter~$k$ has increased by some amount that is divisible
by~$B$. 
When it updates, 
dual agent $q$'s computations take the form
\begin{equation}
\lambda_{[q]}^{q}(tB) = 
\Pi_{M_{q}}[\lambda_{[q]}^{q}(tB-1)+\beta\nabla_{\lambda_{[q]}}L_{\kappa}(z^{q}(tB),\lambda^{q}(tB-1))],
\end{equation}
and between updates it holds the value of~$\lambda^q_q$ constant. 
After updating, dual agent~$q$ sends the new value of~$\lambda^q_q$ to
all primal agents that need it in their computations, and we note
that a primal agent only needs to receive a dual block if
that primal agent's decision variables appear in the constraints
corresponding to those dual variables. Similarly, a dual agent
only needs to receive communications from a primal agent if
that primal agent's decision variables appear in the constraints
encoded in the dual variables that the dual agent updates. 

It is required that all primal agents use the same values
of the dual blocks in their computations. This has been shown
to be necessary~\cite{hendrickson2021totally}
for convergence; any mechanism to
enforce this agreement can be used. 
No communication is required between dual agents since 
updates to~$\lambda^q_q$ do not depend on the
values of other blocks of~$\lambda$. 

We present the full parallelized block gradient-based saddle point algorithm
in Algorithm \ref{alg:algorithm1}. 

\begin{algorithm}[H]
\caption{} \label{alg:algorithm1}

\textbf{Initialize:} All primal variables $z(0)\in Z$, all dual variables
$\lambda(0)\in M$, and the tightening vector $\rho$ as defined in (\ref{eq:RhoEqn}). 
\begin{description}
\item [{for}] $k\in \mathbb{N}$ \textbf{do}
\begin{description}
\item [{for}] $i\in A_{P}$ \textbf{do}
\begin{description}
\item [{if}] $k\in K^{i}$
\begin{description}
\itemsep3pt
\item [{$z_{[i]}^{i}(k+1)=\Pi_{Z_{i}}[z_{[i]}^{i}(k)-\gamma\nabla_{z_{[i]}}L_{\kappa}(z^{i}(k),\lambda^{i}(k))]$}] 

\item [\textbf{Communicate} Agent $i$ sends its updated block to other primal agents, but it may not ] \phantom{a} \newline arrive for some time due to asynchrony.
\end{description}
\item [{else}]~
\begin{description}
\item [{$z_{[i]}^{i}(k+1)=z_{[i]}^{i}(k)$}]~
\end{description}
\item [{end}]~
\item [{if}] $i$ receives $z_{[j]}^{j}$ at time~$k$

$z_{[j]}^{i}(k+1)=z_{[j]}^{j}(\tau_{j}^{i}(k+1))$
\item [{else}]~
\begin{description}
\item [{$z_{[j]}^{i}(k+1)=z_{[j]}^{i}(k)$}]~
\end{description}
\item [{end}]~
\end{description}
\item [{end}]~
\end{description}
\begin{description}
\item [{for}] $q \in A_D$

%
%

\begin{description}
\item [{if}] $k=tB\text{ with }t\in D$
\begin{description}
\item [\textbf{Communicate} All primal agents send their most recent iterate to all dual agents that need it]
\end{description}
\begin{description}
\item [{for}] $q\in A_{D}$ \textbf{do}
\begin{description}
\item [{$\lambda_{[q]}^{q}(tB)=\Pi_{M_{q}}[\lambda_{[q]}^{q}(tB-1)+\beta\nabla_{\lambda_{[q]}}L_{\kappa}(z^{q}(tB),\lambda^{q}(tB-1))]$}]~
\end{description}
\item [{end}]~

\item [\textbf{Communicate} Agent $q$ sends its updated block to the primal agents, but it may not ] \phantom{a} \newline arrive for some time due to asynchrony.
\end{description}
\item [{end}]~
\item [{}]~
\end{description}
\item [{end}]~
\end{description}
\end{description}
\end{algorithm}

\subsection{Convergence Analysis}
We next need to verify that Algorithm \ref{alg:algorithm1} will converge
to the saddle point of $L_{\kappa}$. We start with a proof of convergence
for the primal agents when they have a fixed dual variable onboard.
In it, we define the primal agents' true iterate for all~$k \in \mathbb{N}$ as
\begin{equation}
z(k) = \big(z^1_1(k)^T, z^2_2(k)^T, \ldots, z^m_m(k)^T\big). 
\end{equation}
\begin{lem}
\label{lem:primalConvergence} (Primal Convergence) Let Assumptions
1-3 hold, and consider using Algorithm \ref{alg:algorithm1} to find
a saddle point of $L_{\kappa}$. Fix $t\in D$ and let $t'$ be the
smallest element of D that is greater than t. Fix $\lambda(tB)\in M$
and define $\hat{z}(tB)=\arg\min_{z\in Z}L_{\kappa}(z,\lambda(tB))$
as the point that the primal agents would converge to with $\lambda(tB)$
held fixed. For agents executing Algorithm \ref{alg:algorithm1} at
times $k\in\{tB,tB+1,\ldots,t'B\},$ there exists a scalar $\gamma_{1}>0$
such that for~$\gamma \in (0, \gamma_1)$ we have
\end{lem}
\emph{
\[
\vert\vert z(t'B)-\hat{z}(tB)\vert\vert\leq(1-\theta\gamma)^{t'-t}\vert\vert z(tB)-\hat{z}(tB)\vert\vert,
\]
}

\emph{where $\theta$ is a positive constant, and $(1-\theta\gamma)\in[0,1).$}
\begin{IEEEproof}
We show that Assumptions A and B in \cite{tseng1991rate} are satisfied,
which enables the use of Proposition 2.2 from the same reference to show
convergence.
Assumption A requires that 

1) $L_{\kappa}(\cdot,\lambda(tB))$ is bounded from below on $Z$, 

2) the set $Z$ contains at least one point y such that $y=\Pi_{Z}[y-\nabla_{z}L_{\kappa}(y,\lambda(tB))]$,
and 

3) $\nabla_{z}L_{\kappa}(\cdot,\lambda(tB))$ is Lipschitz on $Z$. 

First, we see that for all~$z \in Z$ we have 
\begin{equation}
\nabla_z L_{\kappa}\big(z, \lambda(tB)\big) = c + \alpha z
+ A^T\lambda(tB),
\end{equation}
which is Lipschitz in~$z$ with constant~$\alpha$. Then Assumption A.3 is satisfied. 


Next, for every choice of $\lambda(tB)\in M$, the function $L_{\kappa}(\cdot,\lambda(tB))$
is bounded from below because it is continuous and its domain~$Z$ is compact. 
Then Assumption A.2 is satisfied. For every $t\in\mathbb{N}$ and
every fixed $\lambda(tB)\in M$, the strong convexity of $L_{\kappa}(\cdot,\lambda(tB))$
implies that is has a unique minimum over $Z$, denoted $\hat{z}(tB)$.
This point is 
the unique fixed point of the projected gradient descent mapping
and Assumption A.1 is satisfied. Then all conditions of Assumption
A in \cite{tseng1991rate} are satisfied. 

Assumption B in \cite{tseng1991rate} requires 

1) the isocost curves of $L_{\kappa}(\cdot,\lambda(tB))$ to be separated,
and 

2) the ``error bound condition'', stated as (2.5) in \cite{tseng1991rate},
is satisfied.

It is observed in~\cite{tseng1991rate} 
that both criteria are satisfied by functions that are strongly convex
with Lipschitz gradients. In this work,  $L_{\kappa}(\cdot,\lambda(tB))$ has
both of these properties. Then Assumption B is satisfied as well, and 
an application of Proposition 2.2 \cite{tseng1991rate} completes the proof.
\end{IEEEproof}
The following lemma is stated here to later be used in the main convergence
theorem for the dual variable. 

\begin{lem}
\label{lem:FromKoshal}Let $\lambda_{1},\lambda_{2}\in\mathbb{R}_{+}^{y}.$
Then for any points $z_{1}$ and $z_{2}$ such that

\begin{equation}
z_{1} = \argmin_{z \in Z} L_{\kappa}(z,\lambda_{1})  
\textnormal{ and }
z_{2} = \argmin_{z \in Z} L_{\kappa}(z,\lambda_{2})
\end{equation}

the pairs $(z_{1},\lambda_{1})$ and $(z_{2},\lambda_{2})$ satisfy

\[
(\lambda_{2}-\lambda_{1})^{T}(-g(z_{2})+g(z_{1}))\geq\frac{\alpha}{\|A\|_2^{2}}\vert\vert g(z_{2})-g(z_{1})\vert\vert^{2}.
\]

Moreover, they satisfy
\begin{equation}
\|\lambda_2 - \lambda_1\| \geq \frac{\alpha}{\|A\|_2}\|z_2 - z_2\|. 
\end{equation}
\end{lem}
\begin{IEEEproof}
See \cite[Lemma 4.1]{koshal2011multiuser}.
\end{IEEEproof}
Next, we derive a proof of convergence for dual agents across two sequential
timesteps.
\begin{lem}
\label{lem:(Dual-Convergence-Between}(Dual Convergence Between Two
Time Steps) Let Assumptions 1-3 hold and consider the use of Algorithm
\ref{alg:algorithm1} to find the saddle point of~$L_{\kappa}$. Let
the dual step size be $0<\beta<\min\begin{Bmatrix}\frac{2\alpha}{\vert\vert A \vert\vert_{2}+2\alpha\delta},\frac{2\delta}{1+\delta^{2}}\end{Bmatrix}$
and let $t_{1}$ and $t_{2}$ denote two consecutive times that dual
updates have occured, with $t_{1}<t_{2}$ and $t_{1},t_{2}\in D.$
Then Algorithm~1 produces dual variables~$\lambda(t_1B)$
and~$\lambda(t_2B)$ that satisfy 
\end{lem}
\begin{equation}
\vert\vert\lambda(t_{2}B)-\hat{\lambda}_{\kappa}\vert\vert^{2}\leq q_{d}\vert\vert\lambda(t_{1}B)-\hat{\lambda}_{\kappa}\vert\vert^{2}+4r^2q_{d}q_{p}^{2(t_{2}-t_{1})}\cdot\vert\vert A\vert\vert_{2}^{2}+8r^2\beta^{2}q_{p}^{t_{2}-t_{1}}\cdot\vert\vert A\vert\vert_{2}^{2}
\end{equation}
\emph{where $q_{d}\triangleq(1-\beta\delta)^{2}+\beta^{2}\in[0,1),\text{ }q_{p}\triangleq(1-\theta\gamma)\in[0,1)$, 
$r$ is from~\eqref{eq:feasProblem-1}, and $\hat{\lambda}_{\kappa}$ is
the dual component of the unique saddle point of $L_{\kappa}.$}
\begin{IEEEproof}
Define $\hat{z}_{\kappa}(t_{1}B)\triangleq\text{arg min}_{z\in Z}L_{\kappa}(z,\lambda(t_{1}B))$
and recall $\hat{z}_{\kappa}=\text{arg min}_{z\in Z}L_{\kappa}(z,\hat{\lambda}_{\kappa})$. 
Given that all dual agents use the same primal variables in
their computations, we analyze all dual agents' computations
simultaneously with the combined dual update law
\begin{equation}
\lambda(t_2B) = \Pi_{M}\big[\lambda(t_1B) + \beta \nabla_{\lambda}L_{\kappa}\big(z(t_2B), \lambda(t_1B)\big)\big],
\end{equation}
where~$z(t_2B)$ is the common primal variable among dual agents. 
Then expanding the dual update law and using the non-expansiveness
of $\Pi_{M}$ gives

\begin{align*}
\vert\vert\lambda(t_{2}B)-\hat{\lambda}_{\kappa}\vert\vert & =\vert\vert\Pi_{M}[\lambda(t_{1}B)+\beta(g(z(t_{2}B)-\delta\lambda(t_{1}B))]-\Pi_{M}[\hat{\lambda}_{\kappa}+\beta(g(\hat{z}_{\kappa})-\delta\hat{\lambda}_{\kappa})]\vert\vert^{2}\\
 &=(1-\beta\delta)^{2}\vert\vert\lambda(t_{1}B)-\hat{\lambda}_{\kappa}\vert\vert^{2}+\beta^{2}\vert\vert g(z(t_{2}B)-g(\hat{z}_{\kappa})\vert\vert^{2}\\
 & -2\beta(1-\beta\delta)(\lambda(t_{1}B)-\hat{\lambda}_{\kappa})^{T}(g(\hat{z}_{\kappa})-g(z(t_{2}B))).
\end{align*}

Adding $g(\hat{z}_{\kappa}(t_{1}B))-g(\hat{z}_{\kappa}(t_{1}B))$
inside the last set of parentheses gives

\begin{align}
\vert\vert\lambda(t_{2}B)-\hat{\lambda}_{\kappa}\vert\vert & = (1-\beta\delta)^{2}\vert\vert\lambda(t_{1}B)-\hat{\lambda}_{\kappa}\vert\vert^{2}+\beta^{2}\vert\vert g(z(t_{2}B)-g(\hat{z}_{\kappa})\vert\vert^{2}\nonumber \\
 & -2\beta(1-\beta\delta)(\lambda(t_{1}B)-\hat{\lambda}_{\kappa})^{T}(g(\hat{z}_{\kappa})-g(\hat{z}_{\kappa}(t_{1}B))\nonumber \\
 & -2\beta(1-\beta\delta)(\lambda(t_{1}B)-\hat{\lambda}_{\kappa})^{T}(g(\hat{z}_{\kappa}(t_{1}B))-g(z(t_{2}B))).\label{eq:Lam_t2_LamOpt}
\end{align}

Applying Lemma \ref{lem:FromKoshal} to the pairs $(\hat{z}_{\kappa},\hat{\lambda}_{\kappa})$
and $(\hat{z}_{\kappa}(t_{1}B),\lambda(t_{1}B))$ results in 

\begin{equation}
(\lambda(t_{1}B)-\hat{\lambda}_{\kappa})^{T}(-g(\hat{z}_{\kappa})+g(\hat{z}_{\kappa}(t_{1}B))) 
 \geq\frac{\alpha}{\vert\vert A\vert\vert_{2}^{2}}\vert\vert g(\hat{z}_{\kappa}(t_{1}B))-g(\hat{z}_{\kappa})\vert\vert^{2},\label{eq:applyingKoshLemma}
\end{equation}
which we apply to the third term on
the right-hand side of \eqref{eq:Lam_t2_LamOpt} to find 

\begin{align}
\vert\vert\lambda(t_{2}B)-\hat{\lambda}_{\kappa}\vert\vert & \leq(1-\beta\delta)^{2}\vert\vert\lambda(t_{1}B)-\hat{\lambda}_{\kappa}\vert\vert^{2}+\beta^{2}\vert\vert g(z(t_{2}B)-g(\hat{z}_{\kappa})\vert\vert^{2}\nonumber \\
 & -2\beta(1-\beta\delta)\frac{\alpha}{\vert\vert A\vert\vert_{2}^{2}}\vert\vert g(\hat{z}_{\kappa}(t_{1}B))-g(\hat{z}_{\kappa})\vert\vert^{2}\nonumber \\
 & -2\beta(1-\beta\delta)(\lambda(t_{1}B)-\hat{\lambda}_{\kappa})^{T}(g(\hat{z}_{\kappa}(t_{1}B))-g(z(t_{2}B))).\label{eq:Lam_t2_LamOpt-1}
\end{align}

To bound the last term, we can see that $0\leq\vert\vert(1-\beta\delta)(g(\hat{z}_{\kappa}(t_{1}B))-g(z(t_{2}B)))+
\beta(\lambda(t_{1}B)-\hat{\lambda}_{\kappa})\vert\vert^{2}$, which
can be expanded and rearranged to give

\begin{align}
-2\beta(1-\beta\delta)(\lambda(t_{1}B)-\hat{\lambda}_{\kappa})^{T}(g(\hat{z}_{\kappa}(t_{1}B))-g(z(t_{2}B))) & \leq(1-\beta\delta)^{2}\vert\vert(g(\hat{z}_{\kappa}(t_{1}B))-g(z(t_{2}B)))\vert\vert^{2}+\beta^{2}\vert\vert(\lambda(t_{1}B)-\hat{\lambda}_{\kappa})\vert\vert^{2}.\label{eq:lastTermBound}
\end{align}

Then we can apply \eqref{eq:lastTermBound} to the
last term in \eqref{eq:Lam_t2_LamOpt-1} to obtain

\begin{align}
\vert\vert\lambda(t_{2}B)-\hat{\lambda}_{\kappa}\vert\vert & \leq(1-\beta\delta)^{2}\vert\vert\lambda(t_{1}B)-\hat{\lambda}_{\kappa}\vert\vert^{2}+\beta^{2}\vert\vert g(z(t_{2}B)-g(\hat{z}_{\kappa})\vert\vert^{2}\nonumber \\
 & -2\beta(1-\beta\delta)\frac{\alpha}{\vert\vert A\vert\vert_{2}^{2}}\vert\vert g(\hat{z}_{\kappa}(t_{1}B))-g(\hat{z}_{\kappa})\vert\vert^{2}\nonumber \\
 & +(1-\beta\delta)^{2}\vert\vert(g(\hat{z}_{\kappa}(t_{1}B))-g(z(t_{2}B)))\vert\vert^{2}+\beta^{2}\vert\vert(\lambda(t_{1}B)-\hat{\lambda}_{\kappa})\vert\vert^{2}.\label{eq:Lam_t2_LamOpt-FinalBound}
\end{align}

Next, add and subtract $g(\hat{z}_{\kappa}(t_{1}B))$ inside the norm
in the second term of \eqref{eq:Lam_t2_LamOpt-FinalBound} to obtain

\begin{align*}
\beta^{2}\vert\vert g(z(t_{2}B))-g(\hat{z}_{\kappa})\vert\vert^{2} & \leq\beta^{2}\vert\vert g(z(t_{2}B))-g(\hat{z}_{\kappa}(t_{1}B))\vert\vert^{2}+\beta^{2}\vert\vert g(\hat{z}_{\kappa}(t_{1}B))-g(\hat{z}_{\kappa})\vert\vert^{2}\\
 & +2\beta^{2}\vert\vert g(z(t_{2}B))-g(\hat{z}_{\kappa}(t_{1}B))\vert\vert\cdot\vert\vert g(\hat{z}_{\kappa}(t_{1}B))-g(\hat{z}_{\kappa})\vert\vert.
\end{align*}

Applying this to the second term in \eqref{eq:Lam_t2_LamOpt-FinalBound}
and grouping terms gives

\begin{align*}
\vert\vert\lambda(t_{2}B)-\hat{\lambda}_{\kappa}\vert\vert & \leq((1-\beta\delta)^{2}+\beta^{2})\vert\vert\lambda(t_{1}B)-\hat{\lambda}_{\kappa}\vert\vert^{2}+\begin{pmatrix}\beta^{2}-2\beta(1-\beta\delta)\frac{\alpha}{\vert\vert A\vert\vert_{2}^{2}}\end{pmatrix}\vert\vert g(\hat{z}_{\kappa}(t_{1}B))-g(\hat{z}_{\kappa})\vert\vert^{2}\\
 & +((1-\beta\delta)^{2}+\beta^{2})\vert\vert g(\hat{z}_{\kappa}(t_{1}B))-g(z(t_{2}B))\vert\vert^{2}\\
 & +2\beta^{2}\vert\vert g(z(t_{2}B)-g(\hat{z}_{\kappa}(t_{1}B))\vert\vert\cdot\vert\vert g(\hat{z}_{\kappa}(t_{1}B))-g(\hat{z}_{\kappa})\vert\vert.
\end{align*}

By hypothesis $0<\beta<\frac{2\alpha}{\vert\vert A\vert\vert_{2}+2\alpha\delta},$ which
makes the second term negative. Dropping this negative term and applying
the Lipschitz property of $g$ give the upper bound

\begin{align*}
\vert\vert\lambda(t_{2}B)-\hat{\lambda}_{\kappa}\vert\vert & \leq((1-\beta\delta)^{2}+\beta^{2})\vert\vert\lambda(t_{1}B)-\hat{\lambda}_{\kappa}\vert\vert^{2}+((1-\beta\delta)^{2}+\beta^{2})\cdot\vert\vert A\vert\vert_{2}^{2}\cdot \vert\vert\hat{z}_{\kappa}(t_{1}B)-z(t_{2}B)\vert\vert^{2}\\
 & +2\beta^{2}\cdot\vert\vert A\vert\vert_{2}^{2}\cdot\vert\vert z(t_{2}B)-\hat{z}_{\kappa}(t_{1}B)\vert\vert\cdot\vert\vert\hat{z}_{\kappa}(t_{1}B)-\hat{z}_{\kappa}\vert\vert\cdot
\end{align*}

Lemma \ref{lem:primalConvergence} can be applied to $\vert\vert\hat{z}_{\kappa}(t_{1}B)-z(t_{2}B)\vert\vert^{2}$
to give an upper bound of
\[
\vert\vert\hat{z}_{\kappa}(t_{1}B)-z(t_{2}B)\vert\vert^{2}\leq q_{p}^{2(t_{2}-t_{1})}\vert\vert z(t_{1}B)-\hat{z}_{\kappa}(t_{1}B)\vert\vert^{2},
\]
where $q_{p}\text{ }=(1-\theta\gamma)\in[0,1)$. Next, the maximum
distance between any two primal variables is bounded by $\underset{z,y\in Z}{\max}\vert\vert z-y\vert\vert \leq \max_{z, y} \|z\| + \|y\| \leq 2r$.
Using this maximum distance and setting $q_{d}=((1-\beta\delta)^{2}+\beta^{2})$,
it follows that $q_{d}\in[0,1)$ because $\beta<\frac{2\delta}{1+\delta^{2}}$.
Finally, we come to the result that

\[
\vert\vert\lambda(t_{2}B)-\hat{\lambda}_{\kappa}\vert\vert^{2}\leq q_{d}\vert\vert\lambda(t_{1}B)-\hat{\lambda}_{\kappa}\vert\vert^{2}+4r^2q_{d}q_{p}^{2(t_{2}-t_{1})}\cdot\vert\vert A\vert\vert_{2}^{2}+8r^2\beta^{2}q_{p}^{t_{2}-t_{1}}\cdot\vert\vert A\vert\vert_{2}^{2}.
\]
\end{IEEEproof}
In brief, the terms that do not explicitly depend on the dual variable
are dependent on regularization parameters and the number of primal agent iterations
that occur between dual agent updates. Now we will use the convergence
rate for two consecutive timesteps and derive a proof of convergence
for the duals to the optimum.
\begin{thm}
\label{thm:dualConvergence} (Dual Convergence to Optimum)\emph{ }Let
Assumptions 1-3 hold and consider the use of Algorithm \ref{alg:algorithm1}
to find the saddle point of~$L_{\kappa}$. Let the dual step size
be $0<\beta<\min\begin{Bmatrix}\frac{2\alpha}{\vert\vert A\vert\vert_{2}+2\alpha\delta},\frac{2\delta}{1+\delta^{2}}\end{Bmatrix}.$\emph{
}Additionally, let $t_{n}$ denote the $n^{\text{th}}$ entry in D
where $t_{1}<t_{2}<...<t_{n}$, which means $t_{n}B$ is the time
when the $n^{\text{th}}$ dual update occurs for all dual agents.
Then the dual agents executing Algorithm \ref{alg:algorithm1} generate
dual variables that satisfy
\end{thm}
\[
 \vert\vert\lambda\left(t_{n}B\right)-\hat{\lambda}_{\kappa}\vert\vert^{2}  \leq q_{d}^{n}\vert\vert\lambda(0)-\hat{\lambda}_{\kappa}\vert\vert^{2}+\left(4r^2q_{d}q_{p}^{2}\cdot \vert\vert A\vert\vert_{2}^{2}+8r^2\beta^{2}q_{p}\cdot \vert\vert A\vert\vert_{2}^{2}\right)\sum_{i=0}^{n}q_{d}^{i},
\]
where $q_{d}=(1-\beta\delta)^{2}+\beta^{2} \in [0, 1),$ $q_{p}=(1-\theta\gamma)\in[0,1)$, $r$ is from~\eqref{eq:feasProblem-1}, 
and $\hat{\lambda}_{\kappa}$ is
the dual component of the unique saddle point of $L_{\kappa}.$
\begin{IEEEproof}
We apply Lemma \ref{lem:(Dual-Convergence-Between} twice to find 
\begin{align*}
\begin{aligned}\vert\vert\lambda\left(t_{n}B\right)-\hat{\lambda}_{\kappa}\vert\vert^{2} & \leq q_{d}\vert\vert\lambda\left(t_{n-1}B\right)-\hat{\lambda}_{\kappa}\vert\vert^{2}+4r^2q_{d}q_{p}^{2\left(t_{n}-t_{n-1}\right)}\cdot \vert\vert A\vert\vert_{2}^{2}+8r^2\beta^{2}q_{p}^{t_{n}-t_{n-1}}\cdot \vert\vert A\vert\vert_{2}^{2}\\
 & \leq q_{d}^{2}\vert\vert\lambda\left(t_{n-2}B\right)-\hat{\lambda}_{\kappa}\vert\vert^{2}+4r^2q_{d}^{2}q_{p}^{2}\cdot \vert\vert A\vert\vert_{2}^{2}+8r^2\beta^{2}q_{d}q_{p}\cdot \vert\vert A\vert\vert_{2}^{2}+4r^2q_{d}q_{p}^{2}\cdot \vert\vert A\vert\vert_{2}^{2}+8r^2\beta^{2}q_{p}\cdot \vert\vert A\vert\vert_{2}^{2},
 \end{aligned}
 \end{align*}
 where the second inequality follows because $t_{n}-t_{n-1}\geq1$.
 Recursively applying Lemma~\ref{lem:(Dual-Convergence-Between} then gives
 \begin{equation}
 \vert\vert\lambda\left(t_{n}B\right)-\hat{\lambda}_{\kappa}\vert\vert^{2}  \leq q_{d}^{n}\vert\vert\lambda(0)-\hat{\lambda}_{\kappa}\vert\vert^{2}+\left(4r^2q_{d}q_{p}^{2}\cdot \vert\vert A\vert\vert_{2}^{2}+8r^2\beta^{2}q_{p}\cdot \vert\vert A\vert\vert_{2}^{2}\right)\sum_{i=0}^{n}q_{d}^{i},
\end{equation}
as desired. 
\end{IEEEproof}
Next we use this bound to formulate a convergence rate for
the primal variables in Algorithm~\ref{alg:algorithm1}. 
\begin{thm}
\label{thm:overallConvergence}(Overall Primal Convergence) Let Assumptions
1-3 hold and the dual stepsize satisfy $0<\beta<\min\{\frac{2\alpha}{\vert\vert A \vert \vert_{2}+2\alpha\delta},\frac{2\delta}{1+\delta^{2}}\}.$
Consider using Algorithm~\ref{alg:algorithm1} to compute the
saddle point of~$L_{\kappa}$, 
and let $t_{n}$ denote
the $n^{\text{th}}$ entry in D where $t_{1}<t_{2}<\ldots<t_{n}$.
That is, $t_{n}B$ is the time at which the $n^{\text{th}}$ dual
update occurs across all dual agents. Then primal agents executing Algorithm
\ref{alg:algorithm1} generate primal iterates that satisfy
\end{thm}
\[
\vert\vert z(t_{n}B)-\hat{z}_{k}\vert\vert\leq 2q_{p}^{t_{n}-t_{n-1}}r +\frac{\vert\vert A\vert\vert_{2}}{\alpha}\left(q_{d}^{n-1}\vert\vert\lambda(0)-\hat{\lambda}_{\kappa}\vert\vert^{2}+\left(4r^2q_{d}q_{p}^{2}\cdot\vert\vert A\vert\vert_{2}^{2}+8r^2\beta^{2}q_{p}\cdot\vert\vert A\vert\vert_{2}^{2}\right)\sum_{i=0}^{n-1}q_{d}^{i}\right)^{1/2},
\]

where $q_{p}=(1-\theta\gamma)\in[0,1),$ $q_{d}$$=(1-\beta\delta)^{2}+\beta^{2}\in[0,1)$\emph{,
and $\hat{\lambda}_{\kappa}$ is the dual at the optimum of $L_{\kappa}.$}
\begin{IEEEproof}
Using the triangle inequality, we have
\begin{align}
\|z(t_nB) - \hat{z}_{\kappa}\| &= \|z(t_nB) - \hat{z}_{\kappa}(t_{n-1}B) + \hat{z}_{\kappa}(t_{n-1}B) - \hat{z}_{\kappa}\| \\
&\leq \|z(t_nB) - \hat{z}_{\kappa}(t_{n-1}B)\|
+ \|\hat{z}_{\kappa}(t_{n-1}B) - \hat{z}_{\kappa}\|. \label{eq:t3tri} 
\end{align}
%
%
%
%
From Lemma~\ref{lem:primalConvergence}, we see that
\begin{align}
\|z(t_nB) - \hat{z}_{\kappa}(t_{n-1}B)\| &\leq (1 - \theta\gamma)^{t_n-t_{n-1}}
\|z(t_{n-1}B) - \hat{z}_{\kappa}(t_{n-1}B)\| \\
  &\leq 2q_p^{t_n - t_{n-1}}r, \label{eq:t3qp}
\end{align}
where the second inequality follows from the definition of~$r$ in~\eqref{eq:feasProblem-1}. 
Next, from Lemma~\ref{lem:FromKoshal}, we see that 
\begin{equation}
\|\hat{z}_{\kappa}(t_{n-1}B) - \hat{z}_{\kappa}\| \leq
\frac{\|A\|_2}{\alpha}\|\lambda(t_{n-1}B) - \hat{\lambda}_{\kappa}\|. \label{eq:t3dualbound}
\end{equation}
Using~\eqref{eq:t3qp} and~\eqref{eq:t3dualbound}
in~\eqref{eq:t3tri} gives
\begin{equation}
\|z(t_nB) - \hat{z}_{\kappa}\| \leq 2q_p^{t_n-t_{n-1}}r
+ \frac{\|A\|_2}{\alpha}\|\lambda(t_{n-1}B) - \hat{\lambda}_{\kappa}\|,
\end{equation}
and bounding the last term using Theorem~\ref{thm:dualConvergence} completes the proof. 
\end{IEEEproof}

\subsection{Mixed-Integer Solution Recovery}

Due to Theorem \ref{thm:SuboptimalityBound}, 
agents computing the saddle point of~$L_{\kappa}$ will compute
a point that is within a bounded distance from the solution to 
Problem \ref{prob:Original-MILP}. 
With the convergence results
in Theorem \ref{thm:overallConvergence}, we have a rate of convergence
to the saddle point and an upper bound on how
far the iterates of Algorithm \ref{alg:algorithm1}
are from the saddle point at any given time. We
round the final output of Algorithm~\ref{alg:algorithm1} to achieve a feasible mixed integer solution,
i.e., a solution that obeys all integrality constraints in Problem~\ref{prob:Original-MILP}. 
Rounding does not drastically effect the feasibility of the solution, and this will
be demonstrated 
in Section \ref{sec:Simulation}. 

An existing error bound on the rounding
of relaxed MILP solutions suggests that a similar theoretical bound is likely to exist
for our saddle point formulation \cite{stein2016error}. 
Given the extensive theoretical developments that such a bound would require,
we defer its development to a future publication. We also observe that, 
for fast computation of MILP solutions, MILP solvers often use heuristics
to speed up exact methods \cite{lodi2013heuristic} and recent applications
have shown effective results from heuristic methods \cite{bertsimas2018multitarget},
\cite{takapoui2020simple}. Therefore, a simple rounding heuristic
is used on the final output generated by Algorithm~\ref{alg:algorithm1} to obtain a
fast solution. One advantage of this rounding heuristic is that it can be executed
onboard each agent individually, without any communication among agents. 
In the next section, we also show that it only introduces slight
suboptimality in the mixed-integer solutions that it produces. 

\section{Simulation}\label{sec:Simulation}

\begin{figure}[H]
\begin{centering}
\includegraphics[scale=0.25]{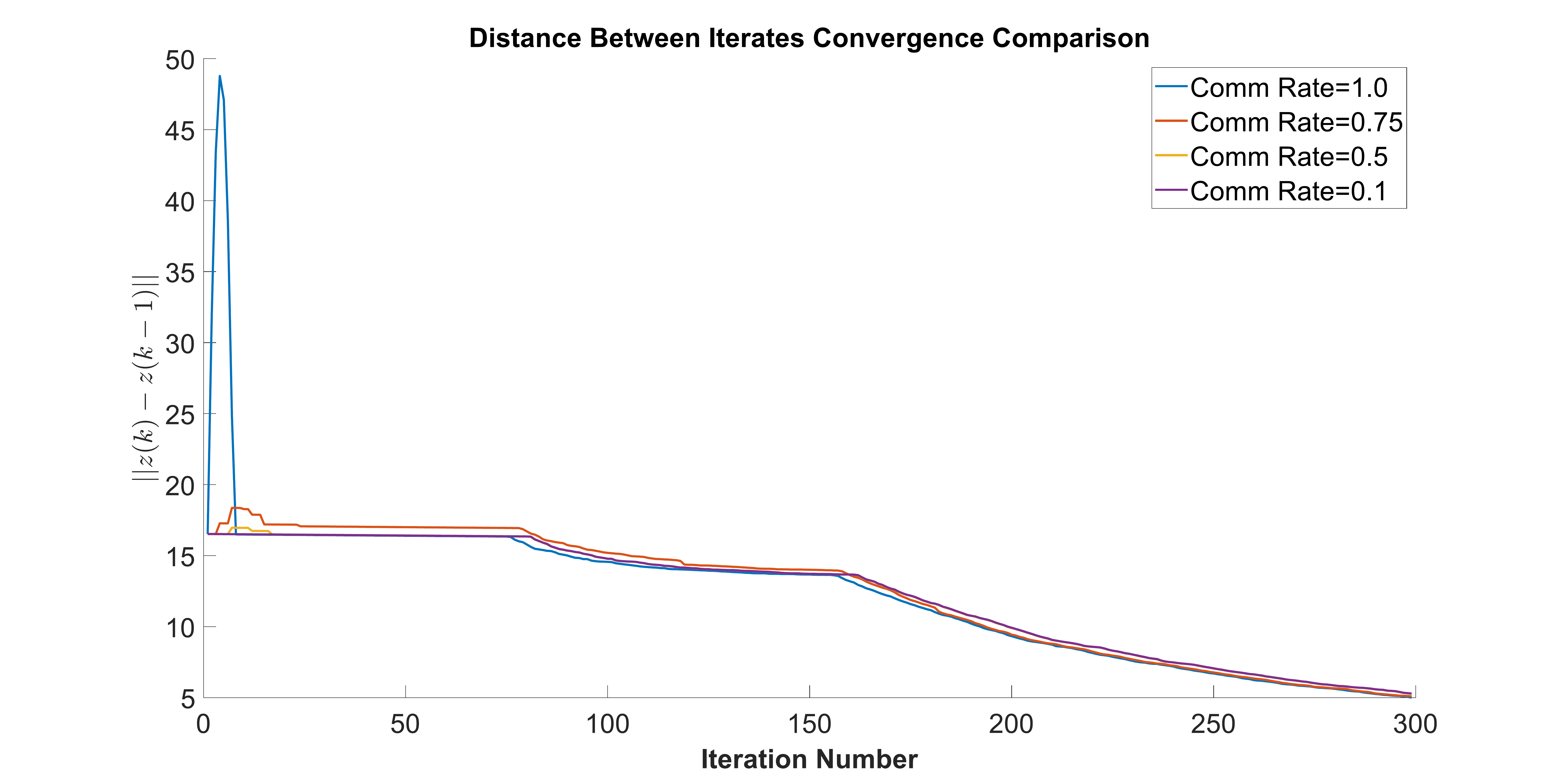}
\par\end{centering}
\centering{}\caption{Four specific representative runs are compared based on communication rates for a random MILP with 300 agents and 30 coupling constraints where there are 285 primal agents and 15 dual agents. The distance between iterates is compared with respect to the communication rates.\label{fig:fig1}}
\end{figure}

\begin{figure}[H]
\begin{centering}
\includegraphics[scale=0.25]{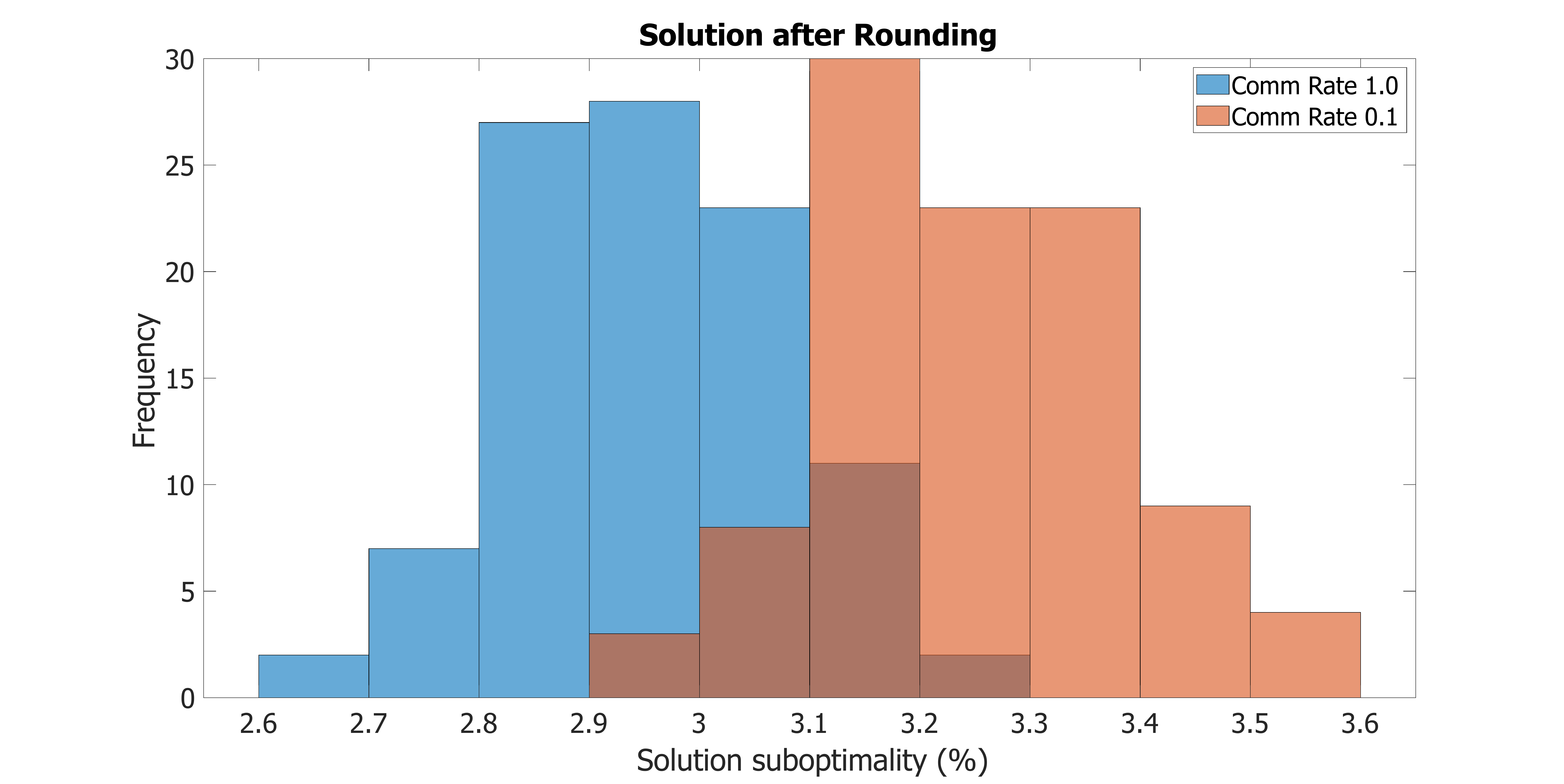}
\par\end{centering}
\centering{}\caption{One hundred
random MILPs were run for 300 agents with 30 coupling constraints where there are 285 primal agents and 15 dual agents in a Monte Carlo simulation. The suboptimality 
in the rounded solutions was computed via~$\frac{\vert c^{T}z^{*}-c^{T}x^{*}|}{\vert c^{T}x^{*}\vert}$.\label{fig:fig2}}
\end{figure}

Random MILPs are used to test Algorithm \ref{alg:algorithm1} in MATLAB
for varying communication rates when computing the saddle point of~$L_{\kappa}$. 
For Algorithm 1, we simulate random communications between agents
using a probability drawn from {[}0,1{]}. The regularization parameters are
$\alpha=10^{-4}$ and $\text{ }\delta=10^{-3}$, and for $300$ agents with 30 coupling
constraints there are 285 primal agents and 15 dual agents. For the
MILPs themselves, the local constraints take the form $-s_{\ell}\leq S_{\ell}x_{\ell}\leq s_{\ell}$, 
where $S_{\ell}=1$ and $s_{\ell}=80$ for all~$\ell \in [285]$. The coupling constraint terms
$A_{\ell}$ for~$\ell \in [285]$ and $b$ are random numbers on the intervals {[}0,1{]} and
{[}20 120{]}, respectively, and the linear cost term $c_{\ell}$ is a random vector on
the interval {[}0 5{]} for all~$\ell \in [285]$. Finally, the local decision vector is $x_{\ell}\in\mathbb{R}^{3}\times\mathbb{Z}^{5}$, 
with a primal step size of $\gamma=0.1.$ One hundred
random MILPs were run in a Monte Carlo simulation, and the suboptimality 
in the rounded solutions was computed via~$\frac{\vert c^{T}z^{*}-c^{T}x^{*}|}{\vert c^{T}x^{*}\vert}$.
For clarity, the convergence of four specific, representative runs is shown in Figure~\ref{fig:fig1},
and the suboptimality of all runs is shown in Figure~\ref{fig:fig2}.

\section{Conclusion\label{sec:Conclusion}}

This paper presents an approximate solution to large-scale mixed-integer
linear programs (MILPs). The MILPs are solved via a distributed saddle
point finding algorithm robust to asynchrony. Theoretical bounds on
the suboptimality gap between the original MILP and the relaxed version
are presented along with convergence rates for the distributed algorithm.
A numerical example verifies the convergence analysis. Future work
will apply the work to large-scale assignment problems and explore
the effect of updating the tightening vector within the algorithm
in an asynchronous setting. 

\bibliographystyle{unsrt}
\bibliography{references}

\end{document}